\documentclass{article}

\title{\LARGE \textbf{Evolutionary \\ Hamiltonian Graph Theory}}
\author{Zh.G. Nikoghosyan\footnote{G.G. Nicoghossian (up to 1997)}}

\begin{document}

\maketitle

\begin{abstract}
We present an alternative domain concerning mathematics to investigate universal evolution mechanisms by focusing on large cycles theory (LCT) - a simplified version of well-known hamiltonian graph theory. LCT joins together a number of $NP$-complete cycle problems in graph theory. $NP$-completeness is the kay factor insuring (by conjecture of Cook) the generation of endless developments and great diversity around large cycles problems. Originated about 60 years ago, the individuals (claims, propositions, lemmas, conjectures, theorems, and so on) in LCT continually evolve and adapt to their environment by an iterative process from primitive beginnings to best possible theorems based on inductive reasoning. LCT evolves much more rapidly than biosphere and has a few thousand pronounced species (theorems). Recall that life on earth with more than 2 million species was originated about 3.7 billion years ago and evolves extremely slowly. We show that all theorems in LCT have descended from some common primitive propositions such as "every complete graph is hamiltonian" or "every graph contains a cycle of length at least one" via improvements, modifications and three kinds of generalizations - closing, associating and extending. It is reasonable to review Darwinian mechanisms in light of LCT evolution mechanisms (especially inductive reasoning) including the origin and macroevolution disputable phenomena in the biosphere.  \\


\end{abstract}

\section{Introduction}

Evolutionary concepts appeared in some early Greek writings, e.g., in the works of Thales, Empedocles, Anaximander, and Aristotle. Under the influence of the Church, no evolutionary theories developed during some 15 centuries. However, much data was accumulated that was to be utilized by later theorists. With the growth of scientific observation and experimentation, there began to appear from about the middle of the 16th century glimpses of the theory of evolution that emerged in the mid 19th century. Charles Robert Darwin set forth the scientific concepts of the evolutionary theory concerning the developments of plants and animals in biosphere that came to be known as Darwinism. In 1859 appeared the first edition of Darwin's "Origin of Species" \cite{[21]}. 

Projections for the total number of species on Earth range from 2 million to 50 million. Scientists have observed that over time, the descendants of living things may change slightly, called "microevolution". Evolutionists teach that small changes accumulated slowly over billions of years and produced the big changes, called "macroevolution".

Darwinism postulates that all organisms on the Earth have descended from a common ancestor over vast periods of time by means of "extremely slight modifications".  To many, this claim sounds reasonable: if small changes can occur within a species, why should not evolution produce big changes over long periods of time. 

Charles Darwin was the first to formulate the theory of evolution by means of natural selection. But how these variations initially arise or are transmitted to offspring, and hence to subsequent generations, was not understood by Darwin.

Natural selection is the only known cause of adaptation, but not the only known cause of evolution. Other, nonadaptable causes of evolution include mutation and genetic drift. In the early 20th century, genetics was integrated with Darwin's theory of evolution by natural selection through the discipline of population genetics. The science of genetics provided a satisfactory explanation for the transmission of variation concerning microevolution. The teaching of macroevolution is built on the claim that mutations - random changes in the genetic code of plants and animals, can produce new species of plants and animals. The gene is the carrier of heredity and determines the attributes of the individual; thus changes in the genes can be transmitted to the offspring and produce new or altered attributes in the new individual. 

One of the hottest controversy in the science is the creation vs. evolution controversy - the Intelligent Design challenge to the theory of evolution. Creationist arguments against evolution theory can sound as follows: 

\begin{itemize}
\item   no one has observed that the small changes we can observe (microevolution) implied that much bigger changes (macroevolution) are also possible,   
\item   scientists worldwide have unearthed and cataloged some 200 million large fossils and billions of small fossils and many researchers agree that this vast and detailed record shows that all the major groups of animals appeared suddenly and remained virtually unchanged, with many species disappearing as suddenly as they arrived.   
\end{itemize}

At the same time, evolutionists argue that creationism is not a scientific theory because it cannot be tested by the scientific method - direct observations and experiments.

This controversy continues to this day and scientists continue to study various aspects of evolution by forming and testing hypotheses, constructing scientific theories, using observational data, and performing experiments in both the field and the laboratory.

Donald Campbell [16] was one of the first authors to formulate a generalized Darwinian concept to explain evolution in a wide variety of other domains, including psychology, economics, culture, language, medicine, computer science and physics.  In result, a number of new areas have been appeared such as evolutionary psychology, culture, language, economics, computation, algorithm, genetic algorithms, programming and so on, inspired mainly by the human intellect. Some of these subjects (psychology, culture, language, economics and so on) trace their origins to the beginning of mankind history and are not well defined. Some others (computation, algorithm, genetic algorithms, programming and so on) are created entirely for practical use to solve optimization problems by direct application of Darwinian mechanisms. 

In this paper we present an alternative domain concerning mathematics to investigate universal evolution mechanisms. We will focus on one of the most heavily studied areas in graph theory that joins together a number of $NP$-complete cycle problems in discrete mathematics, called large cycles theory - a simplified version of well-known hamiltonian graph theory.  

$NP$-completeness is the kay factor insuring (by conjecture of Cook) the generation of endless developments and great diversity around large cycles problems providing an alternative domain comparable with biosphere. We show that the individuals (claims, propositions, lemmas, conjectures, theorems, and so on) in large cycles theory continually evolve and adapt to their environment by an iterative process from simplicity to complexity, from primitive beginnings such as "every complete graph is hamiltonian" or "every graph has a cycle of length at least one" to best possible theorems. We distinguish some evolutionary mechanisms that control this process: improvements, modifications and three kinds of generalizations - closing, associating and extending. Macroevolution (big changes) can be considered as incorporation of a new parameter in a proposition in result of closing or extending generalization in gradual improvement process based on inductive reasoning. Large cycles theory evolves much more rapidly than biosphere and has pronounced species (theorems) with well defined beginning and hereditary microevolution-macroevolution mechanisms. 

As an application, we show that all sharp theorems in large cycles theory have descended from some common ancestors (called fundamental theorems) through extending generalizations. So, large cycles theory provides an exclusive environment where evolution is evident literally, that allows to escape the traditional creation-evolution controversy debates concerning the origin and macroevolution phenomena.  

Large cycles theory can be considered as a simplified evolutionary model concerning human intellect with a number of certain advantages with respect to biology, and its specific evolutionary mechanisms can be useful towards better understanding the evolution mechanisms in biology, as well as the universal mechanisms to explain evolution in a wide variety of domains outside of biology.

\begin{itemize}
\item   Large cycles theory, originated about 60 years ago, has a few thousand pronounced species (theorems) and evolves much more rapidly than living forms on Earth with more than 2 million species, originated about 3.7 billion years ago.
\item   The origins of theorems in large cycles theory can be strongly determined by exact branchings of the tree of improvements and generalizations.
\item   Genetic units and hereditary mechanisms in large cycles theory are much more simpler than gene structures of living forms.
\item   It is quite reasonable to review Darwinian evolutionary mechanisms in light of improvements, modifications and three kinds generalizations - closing, associating and extending, especially inductive reasoning.
\end{itemize}

In the next section, we give necessary terminology and notations. Section 3 is devoted to complexity classes  of computational problems, including $NP$-complete problems. The general environment (large cycles theory) and the patterns evolved in this environment (theorems in large cycles theory) are introduced in Section 4. The structure of theorems and their sharpness are described in Sections 5 and 6.  Evolution mechanisms are classified in Sections 7 and 8. Finally, the definition and the list of "fundamental theorems" are presented in Sections 9 and 10.

\section{Terminology}

Throughout this article we consider only finite undirected graphs without loops or multiple edges. A good reference for any undefined terms is \cite{[14]}.  Let $G$ be a graph with vertex set $V(G)$ and edge set $E(G)$. The degree and the neighborhood of a vertex $x\in V(G)$ are denoted by $d(x)$ and $N(x)$, respectively.

A simple cycle (or just a cycle) $C$ of length $t$ is a sequence $v_1v_2...v_tv_1$ of distinct vertices $v_1,...,v_t$ with $v_iv_{i+1}\in E(G)$ for each $i\in \{1,...,t\}$, where $v_{t+1}=v_1$. When $t=2$, the cycle $C=v_1v_2v_1$ on two vertices $v_1, v_2$ coincides with the edge $v_1v_2$, and when $t=1$, the cycle $C=v_1$ coincides with the vertex $v_1$. So, by this standard definition, all vertices and edges in a graph can be considered as cycles of lengths 1 and 2, respectively.  If  $Q$ is a cycle then we use  $|Q|$ to denote the length of $Q$, that is $|Q|=|V(Q)|$. A path (cycle) on $n$ vertices is denoted by $P_n$ ($C_n$, respectively).

A Hamilton cycle  of a graph is a cycle which passes through every vertex of the graph exactly once, and a graph is hamiltonian if it contains a Hamilton cycle. 

Large cycle structures are centered around well-known Hamilton (spanning) cycles. Other types of large cycles were introduced for different situations when the graph contains no Hamilton cycles or it is difficult to find it. Generally, a cycle $C$ in a graph $G$ is a large cycle if it dominates some certain subgraph structures in $G$ in a sense that every such structure has a vertex in common with $C$. When $C$ dominates all vertices in $G$ then $C$ is a Hamilton cycle. When $C$ dominates all edges in $G$ then $C$ is called a dominating cycle introduced by Nash-Williams \cite{[52]}. Further, if $C$ dominates all paths in $G$ of length at least some fixed integer  $\lambda$ then $C$ is a  $PD_\lambda$(path dominating)-cycle introduced by Bondy \cite{[13]}. Finally, if $C$ dominates all cycles in $G$ of length at least $\lambda$  then $C$ is a  $CD_\lambda$(cycle dominating)-cycle, introduced in \cite{[57]}. 

We reserve $n$, $q$, $\delta$, $\kappa$ and $\alpha$ to denote the number of vertices (order), the number of edges (size), minimum degree, connectivity and the independence number of a graph, respectively. The length $c$ of a longest cycle in a graph is called the circumference. For $C$ a longest cycle in $G$, let $\overline{p}$ and $\overline{c}$ denote the lengths of a longest path and a longest cycle in $G\backslash C$, respectively. 
Put
$$
\delta_t=\min_{v,u}\{\max\{d(v),d(u)\}:v,u\in E(G), d(v,u)=t\}.
$$
Finally, we define $\sigma_t=+\infty$ if $\alpha<\kappa$. Otherwise
$$
\sigma_t=\min\left\{\sum^t_{i=1}d(x_i): \{x_1,x_2,...,x_t\} \ \mbox{is an independent set of vertices in}\ G\right\}.
$$
Let $s(G)$ denote the number of components of a graph $G$. A graph $G$ is $t$-tough if $|S|\ge ts(G\backslash S)$ for every subset $S$ of the vertex set $V(G)$ with $s(G\backslash S)>1$. The toughness of $G$, denoted $\tau(G)$, is the maximum value of $t$ for which $G$ is $t$-tough (taking $\tau(K_n)=\infty$ for all $n\ge 1$). 

Woodall \cite{[66]} defined the binding number $b(G)$ of a graph $G$ as follows:
$$
b(G)=\min_{X\in F}\frac{|N(x)|}{|X|},
$$
where $F=\{X: \emptyset\not=X\subseteq V(G)\}$ and $N(X)=\cup_{x\in X}N(x)$. 

Let $N_{i,j,k}$ be the graph which is obtained by identifying each vertex of a triangle with an endvertex of one of three vertex-disjoint paths of lengths $i,j,k$. If $H_1,...,H_t$ are graphs, then a graph $G$ is said to be $H_1,...,H_t$-free if $G$ contains no copy of any of the graph $H_1,...,H_t$ as induced subgraphs. The graphs $H_1,...,H_t$ will be also referred to in this context as forbidden subgraphs. Denote by $P_t$ the path on $t$ vertices.

A graph $G$ is said to be planar if $G$ is embeddable into the plane without crossing edges. A projective plane, sometimes called a twisted sphere, is a surface without boundary derived from a usual plane by addition of a line at infinity. Just as a straight line in projective geometry contains a single point at infinity at which the endpoints meet, a plane in projective geometry contains a single line at infinity at which the edges of the plane meet. A projective plane can be constructed by gluing both pairs of opposite edges of a rectangle together giving both pairs a half-twist. It is a one-sided surface, but cannot be realized in three-dimensional space without crossing itself.

A graph $G$ is the intersection graph of subgraphs $H_1, . . . ,H_m$ of a graph $H$ if the vertices of $G$ one-to-one correspond
to the subgraphs $H_1, . . . ,H_m$ and two vertices of $G$ are adjacent if and only
if the corresponding subgraphs intersect. 

A graph is an interval graph if and only if it is an intersection graph of subpaths of a path. Next, a graph is a split graph if and only if it is an intersection
graph of subtrees of a star, i.e., a graph $K_{1,m}$. Further, a graph is chordal if and
only if it is an intersection graph of subtrees of a tree. Finally, a comparability graph is a graph whose edges can be transitively oriented
(i.e. if $x > y$ and $y > z$, then $x > z)$; a cocomparability graph $G$ is a graph whose complement $G$
is a comparability graph. Spider graphs are the intersection graphs of subtrees of subdivisions of stars. Thus, spider graphs are chordal graphs that form a common superclass of interval and split graphs.

Let $a,b,t,k$ be integers with $k\le t$. We use $H(a,b,t,k)$ to denote the graph obtained from $tK_a+\overline{K}_t$ by taking any $k$ vertices in subgraph $\overline{K}_t$ and joining each of them to all vertices of $K_b$. Let $L_\delta$ be the graph obtained from $3K_\delta +K_1$ by taking one vertex in each of three copies of $K_\delta$ and joining them each to other. For odd $n\ge 15$, construct the graph $G_n$ from $\overline{K}_{\frac{n-1}{2}}+K_\delta+K_{\frac{n+1}{2}-\delta}$, where $n/3\le\delta\le (n-5)/2$, by joining every vertex in $K_\delta$ to all other vertices and by adding a matching between all vertices in $K_{\frac{n+1}{2}-\delta}$ and $(n+1)/2-\delta$ vertices in $\overline{K}_{\frac{n-1}{2}}$. It is easily seen that $G_n$ is 1-tough but not hamiltonian. A variation of the graph $G_n$, with $K_\delta$ replaced by $\overline{K}_\delta$ and $\delta=(n-5)/2$, will be denoted by $G^*_n$.

\section{Complexity classes of computational problems}

Computational complexity theory focuses on classifying computational problems according to their inherent difficulty, and relating those classes to each other. Graph theory and combinatorics focus on particular problems and their real difficulties. 

Significant progress has been made in combinatorics and graph theory toward improving our understanding of the inherent difficulty in computational problems and what can be computed efficiently. Today, most problems of known interest have been classified as to whether they are polynomial-time solvable or $NP$-complete. 

An algorithm is said to be polynomial time if its running time is upper bounded by a polynomial expression in the size of the input for the algorithm. Problems for which a polynomial time algorithm exists belong to the complexity class $P$, which is central in the field of computational complexity theory. Polynomial time is a synonym for "tractable", "feasible", "efficient", or "fast". The following problems are polynomial-time solvable: shortest path problem, minimum spanning three problem, linear programming, matching, Eulerian cycle problem, network flow problem and so on. 

An algorithm is deterministic if at each step there is only one choice for the next step given the values of the variables at that step. An algorithm is non-deterministic if there is a step that involves parallel processing. A problem is said to be in the class $NP$ of problems if it can be solved by an algorithm which is non-deterministic and has a time complexity function which is polynomial. $NP$ problems are recognized by the fact that their solutions can be checked for correctness by a deterministic polynomial time algorithm. Every problem in $P$ is also in $NP$. The non-deterministic algorithm that can be used is "guess the answer". The guess can be checked in polynomial time by the algorithm which solves the problem. A famous and long standing open problem is whether or not $P = NP$. There is a collection of problems with the property that any polynomial time deterministic algorithm which solves one of them can be converted to a polynomial time algorithm which solves any other one of them (they are said to be polynomially equivalent problems) and if such an algorithm existed for any one of them, then $P = NP$. These problems are called $NP$- hard problems. $NP$-hard problems may or may not be $NP$ problems. Those that are $NP$ are called $NP$-complete problems. An example of an $NP$-complete problem is the Traveling Salesman Problem.

Today, most of important developments in discrete mathematics are centered on various $NP$-complete problems in trying to find different "effective layers" or "effective subspaces" in structures of $NP$-complete problems. In fact, by Cook's conjecture \cite{[19]}, $NP$-complete problems cannot be covered by such layers. Today, after intensive investigations, many $NP$-complete problems are like unbreakable rock fragments with numerous cuttings and bore-holes - tracks of investigations.

\section{General environment and individuals}

Human intellect plays the role of a general environment including large cycles theory as a subarea. Various statements, including claims, propositions, lemmas conjectures and theorems in large cycles theory, play the role of individuals forming a population. 

Large cycles theory traces its origins to 1855. Irish physicist, astronomer and mathematician Sir William Rowan Hamilton (1805-1865)  invented the "Icosian Calculus", a noncommutative algebra so called because it involved a planar embedding of the graph
of a dodecahedron, which has 20 vertices. The system has two operations: $L$ and $R$, standing for "left" and "right" respectively, the idea being that if one has just arrived at a vertex, one can choose to go left or
right, with the value 1 being reserved for an expression which returns to one's point of
origin. For example, a path that turns right twice and then left once can be expressed
as the term $R2L$. Similarly, since each face of a dodecahedron is pentagonal, we know
that $R5 = L5 = 1$. Hamilton showed that symmetry notwithstanding, the equation
$$
LLLRRRLRLRLLLRRRLRLR = 1
$$
defines the only Hamiltonian Cycle on a dodecahedron.
Since $LR\not= RL$, the Icosian Calculus is clearly noncommutative. However, it is
associative. For example, $(LR)L = L(RL)$. Hamilton's first communication about his Icosian Calculus was to
his friend Robert Graves in a letter dated Oct. 7th, 1856.

However, Hamilton   Cycles should not have been
named after Hamilton at all. In fairness, they should be called "Kirkman Cycles" after Thomas Penyngton Kirkman, the man who actually first discovered them. His interest in polyhedra led him to discover Hamilton cycles  in a paper received by the Royal Society on Aug. 6th, 1855, predates
Hamilton's earliest communication, let alone his first publication on the subject, by more
than a year. However, precedence is not the only argument on Kirkman's side. Whereas
Hamilton considered only the one special case of cycles in the dodecahedron, Kirkman's
result was much more general, because he pondered the existence of Hamiltonian Cycles
in all graphs corresponding to planar embeddings of solid shapes. In addition, Kirkman
was the first to discover an infinite class of non-hamiltonian polyhedra. He showed that
any bipartite graph with an odd number of vertices must be non-hamiltonian. He gave
an example of a planar, 3-connected, bipartite, non-hamiltonian graph.

Classic hamiltonian problem; determining when a graph contains a Hamilton cycle, is one of the most central notions in graph theory and is one of the most attractive and most investigated problems among $NP$-complete problems that Karp listed in his seminal paper \cite{[45]}. Cook \cite{[19]} conjectured that one cannot hope for a simple classification of hamiltonian graphs. In other words, it seems to be impossible to obtain a criterion for a graph to be hamiltonian which implies a polynomial-time algorithm. This fact gave rise to a growing number of conditions that are either necessary or sufficient. Today, this conjecture seems much more reasonable motivated by the fact that the developments arising around various $NP$-complete problems in discrete mathematics have undergone a natural gradual growth and evolution, generating a great diversity. This exclusive property of $NP$-complete problems force to think that the diversity arising around such problems potentially should not concede the diversity of living forms in biosphere.  

If a graph $G$ does not satisfy a sufficient condition for hamiltonicity, we cannot guarantee the existence of a Hamilton cycle. But if $G$ is close to satisfy the condition, we may hope find some "hamiltonian-like" structures such as long cycles and hamiltonian paths. Further extensions of these notions lead to cycle and path covers, maximum matching, spanning trees with smallest number of leaves and many others that are rather far from their origins.  Actually, each of these questions is really a part of the general area called "hamiltonian graph theory".

Large cycles theory can be considered as a simplified alternative to hamiltonian graph theory concerning the main "hamiltonian-like" cycle structures in graphs. In fact, large cycles theory is a natural extension of classic hamiltonian problem including Hamilton cycles, longest cycles, dominating cycles, as well as some generalized cycles including Hamilton and dominating cycles as special cases. 

Systematic investigations of Hamilton cycles began only in 1952 when Swiss mathematician Gabriel Andrew Dirac (1925-1984) \cite{[23]} discovered the first sufficient condition for the existence of a Hamilton cycle and the first lower bound for the length of a longest cycle in graphs, based on two simplest graph invariants - order $n$ and minimum degree $\delta$. In the last 60 years, the developments in large cycles theory gave rise to a wide variety of theorems (species, kinds) \cite{[8]}, \cite{[36]}, \cite{[37]}. The following 19 theorems include some generalizations and modifications of Dirac's initial theorems with some progressive tendency, where  $\kappa\le\delta\le\frac{1}{2}\sigma_2\le\frac{1}{3}\sigma_3$ and $\delta\le\delta_2$. \\

\noindent (T1) \ $c\ge\delta+1$   \  \  \  (\cite{[23]}, 1952)

\noindent (T2) \ $\kappa\ge2$  \  \    $\Rightarrow$  \  \    $c\ge\min\{n,2\delta\}$  \  \    (\cite{[23]}, 1952)

\noindent (T3) \ $\kappa\ge2$  \  \    $\Rightarrow$  \  \    $c\ge\min\{n,\sigma_2\}$  \  \    (\cite{[12]}, 1971)

\noindent (T4) \ $\kappa\ge3$, $\delta\ge\alpha$ \  \    $\Rightarrow$  \  \    $c\ge\min\{n,3\delta-3\}$  \  \    (\cite{[39]}, 1978)

\noindent (T5) \ $\kappa\ge3$ \  \    $\Rightarrow$  \  \    $c\ge\min\{n,3\delta-\kappa\}$  \  \    (\cite{[54]}, 1981)

\noindent (T6) \ $\kappa\ge2$  \  \    $\Rightarrow$  \  \    $c\ge\min\{n,2\delta_2\}$  \  \    (\cite{[28]}, 1984)

\noindent (T7) \ $\kappa\ge3$, $G$ is $\delta$-regular \  \    $\Rightarrow$  \  \    $c\ge\min\{n,3\delta\}$  \  \    (\cite{[29]}, 1985) 

\noindent (T8) \ $\kappa\ge 4$, $\delta\ge \alpha$\  \    $\Rightarrow$  \  \   $c\ge\min\{n,4\delta-2\kappa\}$ \ \  (\cite{[55]}, 1985)

\noindent (T9) \ $\tau\ge1$  \  \    $\Rightarrow$  \  \     $c\ge\min\{n,2\delta+2\}$  \  \    (\cite{[4]}, 1986)

\noindent (T10) \ $\tau\ge1$  \  \    $\Rightarrow$  \  \     $c\ge\min\{n,\sigma_2+2\}$  \  \    (\cite{[4]}, 1986)

\noindent (T11) \ $c\ge(\overline{p}+2)(\delta-\overline{p})$ \ \  (\cite{[61]}, 1998)

\noindent (T12) \ $c\ge(\overline{c}+1)(\delta-\overline{c}+1)$ \ \  (\cite{[61]}, 1998)

\noindent (T13) \  $\kappa\ge2$  \  \    $\Rightarrow$  \  \    $c\ge\min\{n,(\tau+1)(\delta+1)-1\}$  \  \    (\cite{[42]}, 1999)

\noindent (T14) \  $\overline{c}\ge\kappa$ \  \    $\Rightarrow$  \  \  $c\ge\frac{(\overline{c}+1)\kappa}{\overline{c}+\kappa+1}(\delta+2)$ \ \  (\cite{[56]}, 2000)

\noindent (T15) \ $\kappa\ge3$ \  \    $\Rightarrow$  \  \    $c\ge\min\{n,\sigma_3-\kappa\}$  \  \    (\cite{[68]}, 2007)

\noindent (T16) \  $\kappa\ge3$, $G$ is claw-free  \  \    $\Rightarrow$  \  \   $c\ge\min\{6\delta-15\}$ \ \  (\cite{[50]}, 2009)

\noindent (T17) \ $\kappa\ge \lambda+2$, $\delta\ge \alpha+\lambda-1$ \      $\Rightarrow$  \     $c\ge\min\{n,(\lambda+2)(\delta-\lambda)\}$ \   (\cite{[57]}, 2009)

\noindent (T18) \ $\kappa\ge 4$, $\delta\ge \alpha$\  \    $\Rightarrow$  \  \   $c\ge\min\{n,4\delta-\kappa-4\}$ \ \  (\cite{[53]}, 2011)

\noindent (T19) \ $\tau>1$ \  \    $\Rightarrow$  \  \     $c\ge\min\{n,2\delta+5\}$  or $G=$Petersen graph \  \     (\cite{[60]}, 2012)\\

As a result of the massive amount of evidence for evolution accumulated in large cycles theory over the last 60 years, we can safely conclude that evolution has occurred and continues to occur in this area.  

Moreover, evolving around an $NP$-complete longest cycle problem, the list of Theorems (T1)-(T19), by conjecture of Cook [19], is not unchanging end-product and will grow  generating continually growing diversity.

\section{The structure of theorems} 

Informally, theorem is of the form of an indicative conditional: 
$$
\mbox{If}\   A  \  \mbox{then} \  B.           \eqno{(1)} 
$$
In this case, $A$ is called the hypothesis (conditions) of Theorem (1) and $B$ the conclusion. Conclusion $B$ indicates the existence of possible types of large cycle structures in a graph $G$. In large cycles theory, conclusion $B$ usually appears in any of the following forms:\\

$(a1)$ $G$ has a Hamilton cycle,

$(a2)$ $G$ has a dominating cycle,

$(a3)$ every longest cycle in $G$ is a dominating cycle,

$(a4)$ $G$ has a $CD_\lambda$-cycle,

$(a5)$ every longest cycle in $G$ is $CD_\lambda$-cycle,

$(a6)$ a lower bound for the circumference.\\

Sometimes, $B\equiv B_1\vee B_2$ where  
$$
B_1\in \{(a1), (a2), (a3), (a4),(a5)\},\ \ B_2\equiv (a6).
$$

As for hypothesis $A$, generally it can be presented as $A_1\wedge A_2 \wedge... \wedge A_m$ where for each $i\in\{1,...,m\}$, $A_i$ appears in the following forms: \\

$(b1)$   $A_i$ is an algebraic (numerical) relation $f_1\ge f_2$ between two algebraic expressions $f_1,f_2$,

$(b2)$   $A_i$ is a structural limitation defined by forbidden subgraphs (examples: forbidden triangle, claw, $P_6$, and so on),

$(b3)$   $A_i$ is a structural limitation defined by direct description (examples: conditions for a graph to be regular, bipartite, interval, chordal, and so on).\\

If $A=A_1\vee A_2$ then Theorem (1) can be partitioned into two independent theorems "if $A_1$ then B" and "if $A_2$ then $B$".

The hypotheses and conclusions defined by $(a_1)-(a6)$ and $(b_1)-(b_3)$, carry the genetic information (genome) of a theorem in forms of initial graph invariants, generalized invariants, forbidden subgraphs and special graph classes. There are a number of well-known basic (initial) invariants of a graph $G$ occurring in various hamiltonian results and having significant impact on large cycle structures, namely order $n$, size $q$, minimum degree  $\delta$, connectivity  $\kappa$, binding number $b(G)$, independence number  $\alpha$, toughness $\tau$  and the lengths of a longest path and a longest cycle in $G\backslash C$ for a given longest cycle $C$, denoted by $\overline{p}$  and  $\overline{c}$, respectively. 

Some of these basic gene elements, especially minimum degree $\delta$, have been generalized (evolved) in terms of degree sequences, degree sums, generalized degree, neighborhood unions and so on, giving rise many generalized theorems.

\section{Relaxation and strengthening} 

Evolutions mechanisms in large cycles theory are based on relaxation and strengthening.\\ 

\noindent \textbf{Definition 1}. Let $f_1\ge f_2$ be a condition in (1) defined by $(b_1)$. We say that the condition $f_1\ge f_2$ can be relaxed in (1) if it can be replaced by $f_1\ge f_2-\epsilon$ for some positive $\epsilon$.\\

\noindent \textbf{Definition 2}. Let "$G$ is $H_1$-free" be a condition in (1) defined by $(b_2)$. We say that "$G$ is $H_1$-free" is stronger than "$G$ is $H_2$-free" if $H_1$ is an induced subgraph of $H_2$.\\

For example, "$G$ is $P_4$-free" is stronger than "$G$ is $P_5$-free" or "$G$ is $P_6$-free". Further,  "$G$ is $N_{0,0,0}$-free" is stronger than "$G$ is $H$-free" for each 
$$
H\in \{N_{0,0,1}, N_{0,0,2}, N_{0,1,1}, N_{0,0,3}, N_{0,1,2}, N_{1,1,1}\}.
$$

If a theorem is not sharp (best possible, tight) then clearly it is incomplete and need further improvement through relaxation and strengthening.\\

\noindent \textbf{Definition 3}. A theorem is said to be sharp in all respects (partly, respectively) if its conclusion cannot be strengthened and each condition (some condition, respectively) in it cannot be relaxed under the same conclusion. \\

According to Definition 3, algebraic relations (see $(b1)$ in previous section) can be gradually (smoothly) relaxed or strengthened forming the best type of hypotheses for relaxing or strengthening.  

Structural limitations defined by forbidden subgraphs (see $(b2)$ in previous section), form the next type of well defined hypotheses in view of relaxing or strengthening. Consider the following theorem based on structural limitations of this type. \\

\noindent\textbf{Theorem A} (Broersma and Veldman \cite{[15]}, 1997).  Every 2-connected $\{K_{1,3},P_6\}$-free graph is hamiltonian.  \\

Generally, it is difficult to check the sharpness related to forbidden subgraphs. However, the following result essentially simplifies this procedure in Theorem A.\\

\noindent\textbf{Theorem B} (Faudree and Gould \cite{[31]}, 1997). Let $R$ and $S$ be connected graphs $(R,S\not=P_3)$ and $G$ be a 2-connected graph of order $n\ge10$. Then $G$ is $(R,S)$-free implies $G$ is hamiltonian if and only if $R=K_{1,3}$ and $S$ is one of the graphs: $P_4, P_5, P_6, N_{0,0,0}, N_{0,0,1}, N_{0,0,2}, N_{0,1,1}, N_{0,0,3}, N_{0,1,2}$ or $N_{1,1,1}$. \\

By Theorem B,  the condition "$G$ is $P_6$-free" in Theorem A cannot be relaxed by replacing it with "$G$ is $H$-free"  for each 
$$
H\in \{P_4, P_5, N_{0,0,0}, N_{0,0,1}, N_{0,0,2}, N_{0,1,1}, N_{0,0,3}, N_{0,1,2}, N_{1,1,1}\}.
$$
Further,  the condition "$G$ is $\{K_{1,3},P_6\}$-free" in Theorem A cannot be relaxed by replacing it with "$G$ is $K_{1,3}$-free" or "$G$ is $P_6$-free" by the following theorem. \\

\noindent\textbf{Theorem C} (Faudree and Gould \cite{[31]}, 1997). Let $R$  be a connected graph  and $G$ be a 2-connected graph. Then $G$ is $R$-free implies $G$ is hamiltonian if and only if $R=P_3$. \\

Finally, the graph $2K_\delta+K_1$ shows that the condition $\kappa\ge2$ in Theorem A  cannot be replaced by $\kappa\ge1$.

So, Theorem A, as well as Theorems 21-25 in Section 9, are best possible.   

Now consider the third type of conditions providing special graph environments (see $(b3)$ in previous section) such as regular, bipartite, interval, chordal, line, spider, split and transitive  graphs, powers of graphs and so on. If the condition cannot be gradually relaxed, it must be removed from the list of conditions as an extraordinary sort of relaxation. Clearly,  $r$-regularity and $(r+1)$-regularity are noncomparable and when we want to relax a condition such as "$G$ is $r$-regular", we have to remove this condition.

By relaxing the condition "$G$ is bipartite" we get a trivial case when "$G$ is one-partite" or empty graph.

Planarity can be interpreted both in view of forbidden subgraphs and embedding in a plane without crossings. The following well-known theorem is similar to Theorem G and shows that in both cases we get a non planar graph when we try to relax the planarity condition .\\

\noindent \textbf{Theorem D} (Kuratowski \cite{[49]}, 1930). A graph is planar if and only if it does not contain a subgraph that is homeomorphic to $K_5$ or $K_{3,3}$.

\section{Evolution mechanisms in large cycles theory}	 

All theorems in large cycles theory have descended from trivial (primitive) propositions such as:\\

$(c1)$   every complete graph is hamiltonian,

$(c2)$   every graph contains a cycle of length at least one,

$(c3)$   every graph with $n=1,2,...,10$ and $\delta\ge n-1$ is hamiltonian,

$(c4)$   every graph with $\delta\ge n-1\ge1$ is hamiltonian,

$(c5)$   every graph with $\alpha\le1$ is hamiltonian, 

$(c6)$   every graph with $q\ge n(n-1)/2$ is hamiltonian,\\

\noindent via the following evolutionary mechanisms:\\

\begin{itemize}
\item   improvements (vertical evolution),
\item   modifications (horizontal evolution),
\item   closing generalizations,
\item   associating generalizations,
\item   extending generalizations.
\end{itemize}

\subsection{Improvements}	 

Improvement is a progressive (vertical) iterative process in evolution toward finding better results. \\

\noindent \textbf{Definition 4}. Improvement is one of the following procedures:

$(d1)$   relaxing one of the conditions and preserving the conclusion, 

$(d2)$   strengthening the conclusion and preserving the conditions. \\

Improvements are applicable only to the trivial or incomplete (not sharp) results such as $(c3)-(c6)$. For example,  $(c3)$ can be iteratively improved to "every graph with $n=10$ and $\delta\ge i$ is hamiltonian" for $i=8,7,6,5$. The best result in this process can be formulated as follows\\

$(c7)$   every graph with $n=10$ and $\delta\ge5$ is hamiltonian.\\

Furthermore, $(c7)$ can be iteratively improved to "every graph with $n=i$ and $\delta\ge i/2$ is hamiltonian" for $i=11,12,...$.

\subsection{Modifications}	 

Modification is a horizontal developmental process in evolution generating noncomparable results.  \\

\noindent \textbf{Definition 5}. Modification is one of the following procedures:

$(e1)$    relaxing of some conditions, at the same time strengthening some others, under the same conclusion,

$(e2)$    relaxing of some conditions, at the same time relaxing the conclusion,

$(e3)$    strengthening of some conditions, at the same time strengthening the conclusion.  \\

Observing that $\tau\ge 1$ is stronger than $\kappa\ge 2$, and "$G$ is hamiltonian" is stronger than "$G$ contains a dominating cycle", we can state that the following theorems, by Definition 5, are modifications.\\

$(f1)$   If $\kappa\ge 2$ and $\delta\ge (n+2)/3$ then $G$ contains a dominating cycle.   \cite{[52]}

$(f2)$   If $\tau\ge 1$ and $\delta\ge n/3$ then $G$ contains a dominating cycle.   \cite{[9]}

$(f3)$   If $\kappa\ge 2$ and $\delta\ge (n+\kappa)/3$ then $G$ is hamiltonian.   \cite{[54]}\\

\subsection{Closing generalizations}	 

\noindent \textbf{Definition 6}. Closing generalization is an improvement process which yields a best possible result. Often, it is based on inductive reasoning which generates new parameters. \\

As noted above, $(c7)$ can be iteratively improved to "every graph with $n=i$ and $\delta\ge i/2$ is hamiltonian" for $i=11,12,...$. Inductive reasoning allows to obtain a best possible theorem involving the order $n$ as a new parameter. \\

$(c8)$   If $\delta\ge n/2$ then $G$ is hamiltonian.   \cite{[23]}\\

Inductive reasoning is also known as induction: a kind of reasoning that constructs or evaluates propositions that are abstractions of observations of individual instances.

\subsection{Associating generalizations}	 

\noindent \textbf{Definition 7}. Associating generalization joins together closely related noncomparable results for special values of some parameter $\lambda=1,2,...$, based on inductive reasoning. \\

For example, the following theorem \\

$(g1)$   if $\kappa\ge \lambda\ge1$  and $\delta\ge (n+2)/(\lambda+1)+\lambda-2$ then $G$ contains a  $CD_{\min\{\lambda,\delta-\lambda+1\}}$-cycle  \cite{[57]},\\

\noindent associates the following noncomparable results for $\lambda=1,2,3$.\\

$(g2)$   if $\kappa\ge 1$ and $\delta\ge n/2$ then $G$ is hamiltonian \cite{[23]},

$(g3)$   if $\kappa\ge 2$ and $\delta\ge (n+2)/3$  then $G$ contains a dominating cycle \cite{[52]},

$(g4)$   if $\kappa\ge 3$ and $\delta\ge (n+6)/4$ then $G$ contains a $CD_3$-cycle \cite{[41]}.\\

\subsection{Extending generalizations}	 

Since the best possible theorems cannot be improved, they can be involved only by extensions of some notions.\\

\noindent \textbf{Definition 8}. Extensions of some concepts in best possible theorems generate a new kind of so called extending generalizations.\\

The concepts $\sigma_t$ and $\delta_t$ $(t\ge1)$ are two extensions of the minimum degree $\delta$ with $\sigma_1=\delta_1=\delta$ and therefore, the following two theorems \\

$(h1)$   if $\sigma_2\ge n$ then $G$ is hamiltonian  \cite{[63]}, 

$(h2)$   if $\delta_2\ge n/2$ then $G$ is hamiltonian  \cite{[28]}, \\

\noindent are two extending generalizations of $(c8)$.

\section{Microevolution and macroevolution}

Microevolution in large cycles theory is a gradual improvement process based on numerical expressions. 

Macroevolution occurs in result of a closing generalization which determines a new result with a new involved parameter based on inductive reasoning. In other words, macroevolution is a transition from numerical subexpressions to parametrical subexpressions.

All intermediate microevolution changes in improvement process are immediately forgotten and eliminated preserving only the final result as a macroevolution big change. In result, the evolution process in large cycles theory seems discrete with large breaking-offs. Conversely, each new parameter can be incorporated into the theorem a result of some closing generalization. Observe that along with new involved parameters, improvement process preserves the hereditary information in forms of earlier involved parameters.

Graph invariants and their various extensions, combined in convenient relations as parameters, contain global and general information about a graph and its cycle structures like gene structures. They form the hereditary information of theorems in large cycles theory. 

The order $n$ and size $q$ as gene elements one by one are neutral graph invariants with respect to cycle structures. Meanwhile, they become more active combined together (as in Theorem 1). 

The minimum degree  $\delta$ plays a central role in majority of hamiltonian results. It is not too primitive and not too complicated, becoming the most flexible invariant for various possible generalizations. Minimum degree is a more essential invariant than the order and size, providing some dispersion of the edges in a graph. The combinations between order n and minimum degree   become much more fruitful especially under some additional connectivity conditions. 

The impact of some relations on cycle structures can be strengthened under additional conditions of the type $\delta\ge\alpha \pm i$  if for appropriate integer  $i$. Determining the independence number  $\alpha$ is shown in \cite{[35]} to be $NP$-hard problem.

Connectivity is the most valuable research tool toward cognation of large cycle structures. In \cite{[25]}, it was proved that connectivity  $\kappa$ can be determined in polynomial time. Many graph theorists think that the connectivity   is at the heart of all path and cycle questions providing comparatively more uniform dispersion of the edges. 

The binding number $b(G)$ is a measure of how well-knot a graph is. Like the connectivity, the binding number also can be computed in polynomial time, using network techniques \cite{[20]}. 

An alternate connectedness measure is toughness $\tau$  - the most powerful and less investigated graph invariant introduced by Chv\'{a}tal \cite{[17]} as a means of studying the cycle structure of graphs. Moreover, it was proved \cite{[2]} that for any positive rational number $t$, recognizing $t$-tough graphs (in particular 1-tough graphs) is an $NP$-hard problem. Chv\'{a}tal \cite{[17]} conjectured that there exists a finite constant $i_0$  such that every  $i_0$-tough graph is hamiltonian. This conjecture is still open. 

For a given cycle $C$, the idea of using $G\backslash C$ appropriate structures lies in the base of almost all existing proof techniques in trying to construct longer cycles in graphs by the following standard procedure: choose an initial cycle  $C_0$  in $G$ and try to enlarge it by replacing a segment $P^\prime$  of $C_0$  with a suitable path  $P^{\prime\prime}$ longer than  $P^\prime$, having the same end vertices and passing through  $G\backslash C_0$. To find suitable $P^\prime$  and  $P^{\prime\prime}$, one can use the paths or cycles (preferably large) in  $G\backslash C_0$ and connections (preferably high) between these paths (cycles) and  $C_0$. The latter are closely related to  $\overline{p}, \overline{c}$, as well as minimum degree $\delta$  (local connections) and connectivity $\kappa$  (global connections). 

Forbidden small subgraphs provide the next powerful gene element of structural nature that directly force the graph to have large cycles. For example, 2-connected $P_3$-free graphs are hamiltonian since they are complete graphs. The most common of forbidden subgraphs is the claw $K_{1,3}$.

Finally, some special graph classes, that can be defined by direct description,  provide convenient environments to construct large cycles in graphs. They are  regular graphs, planar graphs, bipartite graphs, chordal graphs, interval graphs and so on.

\section{On fundamental results in large cycles theory}

What makes a theorem (problem, conjecture) beautiful?  By G.H. Hardy, "The mathematician's patterns, like the painter's or the poet's
must be beautiful; the ideas, like the colors or the words must fit
together in a harmonious way. Beauty is the first test: there is
no permanent place in this world for ugly mathematics".

In \cite{[11]}, Bondy introduced some criteria to classify conjectures, which can be applicable for theorems as well:

\begin{itemize}
\item   Simplicity: short, easily understandable statement relating
basic concepts.
\item   Element of Surprise: links together seemingly disparate
concepts.
\item   Generality: valid for a wide variety of objects.
\item   Centrality: close ties with a number of existing theorems
and/or conjectures.
\item   Longevity: at least twenty years old.
\item   Fecundity: attempts to prove the conjecture have led to new
concepts or new proof techniques.
\end{itemize}

However, the first formal criterion toward classifying the theorems and conjectures is the property to be best possible (sharp, tight), widely applicable in combinatorics and graph theory. This criterion after some improvement can be applicable in other areas of science.
 
The next formal criterion to distinguish some special kind of theorems is presented in \cite{[62]} by focusing on pure relations between  simplest graph invariants and large cycles structures. These simplest kind of relations having no forerunners in the area, actually form a source from which nearly all possible hamiltonian results can be developed further by various additional new ideas, generalizations, extensions, restrictions and structural limitations.

In this paper we introduce the third formal criterion to distinguish some top theorems in large cycles theory called "fundamental" based on all exact branchings of the tree of generalizations. By this approach, all results in large cycles theory have descended from a number of common ancestors (fundamental result) through extending generalizations. Fundamental results cannot be directly improved and can be evolved only by modifications and generalizations (associating and extending). 

The term "fundamental result" is used in various fields of science to characterize mainly the central and most important results in the area, based on subjective perception. In this paper, this term is used according to the second much more important mean: "forming the source or base from which everything else is made; not able to be divided any further". Observe also that in general, there are no physical and abstract units in the nature, lying in the base of all material or abstract notions. However, every notion in large cycles theory has certain origins due to certain frames of this theory.

\section{The list of fundamental results}

\subsection{Hamilton cycles}

\noindent\textbf{Theorem 1} (Erd\"{o}s and Gallai, 1959) \cite{[27]} 

\noindent Every graph is hamiltonian if
$$
q\ge\frac{n^2-3n+5}{2}.
$$

Example for sharpness. To see that the size bound $(n^2-3n+5)/2$ in Theorem 1 is best possible,  note that the graph formed by joining one vertex of $K_{n-1}$ to $K_1$, contains $(n^2-3n+4)/2$ edges and is not hamiltonian. \\

\noindent\textbf{Theorem 2} (Erd\"{o}s, 1962) \cite{[26]} 

\noindent Every graph is hamiltonian if $1\le\delta\le n/2$ and 
$$
q>\max\left\{\frac{(n-\delta)(n-\delta-1)}{2}+\delta^2,\frac{\left(n-\lfloor\frac{n-1}{2}\rfloor\right)\left(n-\lfloor\frac{n-1}{2}\rfloor-1\right)}{2}+\left\lfloor\frac{n-1}{2}\right\rfloor^2\right\}.
$$

Example for sharpness. The graph consisting of a complete graph on $n-\delta$ vertices, $\delta$ of which are joined to each of $\delta$ independent vertices, shows that the condition in Theorem 2 cannot be weakened.\\

\noindent\textbf{Theorem 3} (Moon and Moser, 1963) \cite{[51]}

\noindent Every balanced bipartite graph is hamiltonian if 
$$
q\ge\frac{n^2-2n+5}{4}.
$$

Examples for sharpness. Clearly, the condition "$G$ is balanced" in Theorem 3 cannot be removed. The graph obtained from $K_{t,t}$ by deleting $t-1$ edges with a common vertex, shows that  the condition $q\ge(n^2-2n+5)/4$ in Theorem 3 cannot be replaced by $q\ge(n^2-2n+4)/4$. \\

\noindent\textbf{Theorem 4} (Moon and Moser, 1963) \cite{[51]}

\noindent Every balanced bipartite graph is hamiltonian if 
$$
q > \frac{n(n-2\delta)}{4}+\delta^2.
$$

Examples for sharpness. Clearly, the condition "$G$ is balanced" in Theorem 4 cannot be removed. Consider the balanced bipartite graph $G=(X,Y;E)$ with vertex classes of the form $X=P\cup Q$, $Y=R\cup S$, where $|P|=|R|=\delta$, $|Q|=|S|=n/2-\delta$, $N_G(x)=R$ for all $x\in P$, and $N_G(x)=Y$ for all $x\in Q$. This example shows that Theorem 4 is best possible.\\

\noindent\textbf{Theorem 5} (Nikoghosyan, 2011) \cite{[58]} 

\noindent Every graph is hamiltonian if
$$
q\le \delta^2+\delta-1.
$$

Example for sharpness. $K_1+2K_{\delta}$.\\

\noindent\textbf{Theorem 6} (Dirac, 1952) \cite{[23]} 

\noindent Every graph is hamiltonian if
$$
\delta\ge\frac{n}{2}.
$$

Example for sharpness. $2K_\delta+K_1$. \\

\noindent\textbf{Theorem 7} (Moon and Moser, 1963) \cite{[51]}

\noindent Every balanced bipartite graph is hamiltonian if
$$
\delta\ge\frac{n+1}{4}.
$$

Examples for sharpness. Clearly, the condition "$G$ is balanced" in Theorem 7 cannot be removed. Since $n$ is even,  the condition $\delta\ge(n+1)/4$ in Theorem 7 yields a stronger condition $\delta\ge(n+2)/4$.  Let $P_i=x_iy_iz_iw_i$ $(i=1,2,3)$ be three disjoint paths. Form a graph from $P_1,P_2,P_3$ by identifying $x_1,x_2,x_3$ in one vertex and $w_1,w_2,w_3$ in another vertex. The resulting graph shows that the condition $\delta\ge(n+1)/4$ in Theorem 7 cannot be replaced by $\delta\ge n/4$. \\

\noindent\textbf{Theorem 8} (Jung, 1978) \cite{[39]} 

\noindent Every graph is hamiltonian if $n\ge 11$, $\tau\ge 1$ and 
$$
\delta\ge\frac{n-4}{2}.
$$

Examples for sharpness. Petersen graph; $K_{\delta, \delta+1}$; $G^*_n$.\\

\noindent\textbf{Theorem 9} (Nikoghosyan, 2012) \cite{[60]} 

\noindent Every graph is hamiltonian if  $\tau > 4/3$ and 
$$
\delta\ge\frac{n-5}{2}.
$$

Examples for sharpness. The Petersen graph shows that the condition $\tau > 4/3$ in Theorem 9 cannot be replaced by $\tau = 4/3$. Let $H_1$ be a complete graph with vertex set $V(H_1)=\{x_1,x_2,x_3,x_4,x_5\}$ and  $H_2$ a complete bipartite graph with bipartition $(V_1,V_2)$, where $V_1=\{y_1,y_2,y_3,y_4,y_5\}$ and $|V_2|=2$.   The graph obtained from disjoint graphs $H_1$ and $H_2$ by adding the edges $x_iy_i$ $(i=1,2,3,4,5)$, shows that the condition $\delta\ge(n-5)/2$ in Theorem 9 cannot be replaced by $\delta\ge(n-6)/2$.\\

\noindent\textbf{Theorem 10} (Nikoghosyan, 1981) \cite{[54]} 

\noindent Every graph is hamiltonian if $\kappa \ge 2$ and 
$$
\delta\ge \frac{n+\kappa}{3}.
$$

Examples for sharpness. $2K_\delta+K_1$; $H(1,\delta-\kappa+1,\delta,\kappa)$ $(2\le\kappa<n/2)$.\\

\noindent\textbf{Theorem 11} (Bauer and Schmeichel, 1991) \cite{[5]} 

\noindent Every graph is hamiltonian if $\tau\ge 1$ and 
$$
\delta\ge \frac{n+\kappa-2}{3}.
$$

Examples for sharpness. $K_{\delta,\delta+1}$; $L_\delta$.\\

\noindent\textbf{Theorem 12} (Nash-Williams, 1971) \cite{[52]} 

\noindent Every graph is hamiltonian if $\kappa\ge 2$ and 
$$
\delta\ge\max\left\{\frac{n+2}{3},\alpha\right\}.
$$

Examples for sharpness. $(\lambda+1)K_{\delta-\lambda+1}+K_\lambda$ $(\delta\ge2\lambda)$; $(\lambda+2)K_{\delta-\lambda}+K_{\lambda+1}$ $(\delta\ge2\lambda+1)$; $H(\lambda,\lambda+1,\lambda+3,\lambda+2)$.\\

\noindent\textbf{Theorem 13} (Bigalke and Jung, 1979) \cite{[9]}

\noindent Every graph is hamiltonian if $\tau\ge 1$ and 
$$
\delta\ge\max \left\{\frac{n}{3},\alpha-1\right\}.
$$

Examples for sharpness. $K_{\delta,\delta+1}$ $(n\ge 3)$; $L_\delta$ $(n\ge 7)$; $K_{\delta,\delta+1}$ $(n\ge 3)$.\\

\noindent\textbf{Theorem 14} (Fraisse, 1986) \cite{[34]} 

\noindent Let $G$ be a graph and $\lambda$ a positive integer. Then $G$ is hamiltonian if  $\kappa\ge\lambda+1$ and 
$$
\delta\ge\max\left\{\frac{n+2}{\lambda+2}+\lambda-1,\alpha+\lambda-1\right\}.
$$

Examples for sharpness. $(\lambda+1)K_{\delta-\lambda+1}+K_\lambda$ $(\delta\ge2\lambda)$; $(\lambda+2)K_{\delta-\lambda}+K_{\lambda+1}$ $(\delta\ge2\lambda+1)$; $H(\lambda,\lambda+1,\lambda+3,\lambda+2)$. Theorem 14 can be considered as a union (not a generalization) of fundamental results for all possible values of $\lambda$.\\

\noindent\textbf{Theorem 15} (Yamashita, 2008) \cite{[69]}

\noindent Every graph is hamiltonian if $\kappa\ge 3$ and 
$$
\delta\ge \max\left\{\frac{n+\kappa+3}{4}, \alpha\right\}.
$$

Examples for sharpness. $3K_{\delta-1}+K_2$; $H(2,n-3\delta+3,\delta-1,\kappa)$; $H(1,2,\kappa+1,\kappa)$.\\

\noindent\textbf{Theorem 16} (Chv\'{a}tal and Erd\"{o}s, 1972) \cite{[18]} 

\noindent Every graph is hamiltonian if
$$
\kappa\ge \alpha.
$$

Example for sharpness. $K_{\delta,\delta+1}$.\\

\noindent\textbf{Theorem 17} (Woodall, 1973) \cite{[66]} 

\noindent Every graph $G$ is hamiltonian if
$$
b(G)\ge \frac{3}{2}.
$$

Example for sharpness. $aK_2+\overline{K}_{a-1}$.\\

\noindent\textbf{Theorem 18} (Fleischner, 1974) \cite{[33]}

\noindent The square of every 2-connected graph is hamiltonian.  \\

Examples for sharpness. Clearly, the power of a graph cannot be reduced to one in Theorem 18, since there are 2-connected nonhamiltonian graphs. Next, 2-connectivity condition in Theorem 18 cannot be relaxed  since the square of a graph $G$ is not hamiltonian if  $G-x$ has at least three nontrivial components in which $x$ has exactly one neighbor.\\

\noindent\textbf{Theorem 19} (Tutte, 1956) \cite{[65]}  

\noindent Every  4-connected planar graph is hamiltonian.\\

Examples for sharpness. Tutte's graph shows that 4-connectivity condition in Theorem 19 cannot be relaxed. Complete bipartite graph $K_{4,5}$ shows that planarity is a necessary condition in Theorem 19.\\

\noindent\textbf{Theorem 20} (R. Thomas and X. Yu, 1994) \cite{[64]}  

\noindent Every  4-connected projective-plane graph is hamiltonian.\\

Examples for sharpness. The simplest non-orientable surface on which the Petersen graph can be embedded without crossings is the projective plane. The Petersen graph shows that 4-connectivity condition in Theorem 20 cannot be relaxed. On the other hand, there are 4-connected non hamiltonian graphs that cannot be embedded on projective plane (otherwise, all 4-connected graphs are hamiltonian), implying that the condition "$G$ is projective plane graph" cannot be removed in Theorem 20. \\

\noindent\textbf{Theorem 21} (Faudree and Gould, 1997) \cite{[31]}

\noindent Every 2-connected $P_3-$free graph is hamiltonian.  \\

Examples for sharpness. See Subsection 5.3.\\

\noindent\textbf{Theorem 22} (Broersma, Veldman, 1997) \cite{[15]}

\noindent Every 2-connected $\{K_{1,3},P_6\}-$free graph is hamiltonian.  \\

Examples for sharpness. See Subsection 5.3.\\

\noindent\textbf{Theorem 23} (Faudree, Gould, Ryj\'{a}\v{c}ek and Schiermeyer, 1997) \cite{[32]}

\noindent Every 2-connected $\{K_{1,3},N_{0,0,3}\}-$free graph with $n\ge 10$ is hamiltonian.  \\

Examples for sharpness. See Subsection 5.3.\\

\noindent\textbf{Theorem 24} (Bedrossian, 1997) \cite{[7]}

\noindent Every 2-connected $\{K_{1,3},N_{0,1,2}\}-$free graph is hamiltonian.  \\

Examples for sharpness. See Subsection 5.3.\\

\noindent\textbf{Theorem 25} (Duffus, Jakobson and Gould, 1997) \cite{[24]}

\noindent Every 2-connected $\{K_{1,3},N_{1,1,1}\}-$free graph is hamiltonian.  \\

Examples for sharpness. See Subsection 5.3.\\

\noindent\textbf{Theorem 26} (Keil, 1985) \cite{[46]}

\noindent Every 1-tough interval graph is hamiltonian.\\

Examples for sharpness. Star graphs are interval nonhamiltonian graphs with $\tau<1$, implying that 1-toughness condition in Theorem 26 cannot be relaxed. The Petersen graph shows that the condition "$G$ is interval graph" in Theorem 26 cannot be removed. \\

\noindent\textbf{Theorem 27}  (Kratsch, Lehel and M\"{u}ller, 1996) \cite{[48]}

\noindent Every $3/2-$tough split graph is hamiltonian.\\

Examples for sharpness. In \cite{[48]}, $(3/2-\epsilon)$-tough split graphs are constructed that are not hamiltonian. There are non hamiltonian graphs with $\tau=9/4-\epsilon>3/2$, implying that the condition "$G$ is split graph" in Theorem 27 cannot be removed.\\

\noindent\textbf{Theorem 28} (Deogun, Kratsch and Steiner, 1997) \cite{[22]}

\noindent Every 1-tough cocomparability graph is hamiltonian.\\

Examples for sharpness. Clearly, any complete graph is a comparability graph and hence, any empty graph is a cocomparability graph with $\tau<1$, implying that the condition "$G$ is 1-tough" in Theorem 28  cannot be relaxed. On the other hand, there are 1-tough non hamiltonian non cocomparability graphs (otherwise, all 1-tough graphs are hamiltonian), implying that the condition "$G$ is cocomparability graph" in Theorem 28 cannot be removed.\\

\noindent\textbf{Theorem 29}  (B\"{o}hme, Harant and Tk\'{a}\v{c}, 1999) \cite{[10]}

\noindent Every chordal, planar graph with $\tau>1$ is hamiltonian. \\

Examples for sharpness. In \cite{[10]}, it is proved that for any $\epsilon > 0$, there is a 1-tough chordal planar graph $G_\epsilon$ such that the length of a longest cycle of $G_\epsilon$ is less than $\epsilon|V(G_\epsilon)|$, implying that the condition $\tau>1$ in Theorem 29 cannot be relaxed. Chv\'{a}tal \cite{[17]} obtained $(3/2-\epsilon)$-tough graphs without a 2-factor, implying that the planarity condition in Theorem 29 cannot be removed. Finally, Harant \cite{[38]} found 3/2-tough planar nonhamiltonian graphs, implying that the condition "$G$ is chordal" in Theorem 29 cannot be removed.   \\

\noindent\textbf{Theorem 30} (Kaiser, Kr\'{a}l and Stacho, 2007) \cite{[43]}

\noindent Every $3/2$-tough spider (intersection) graph is hamiltonian.  \\

Examples for sharpness. In \cite{[43]}, Kaiser, Kr\'{a}l and Stacho constructed $(3/2-\epsilon)$-tough spider graphs that do not contain a Hamilton cycle, implying that the condition "$G$ is 3/2-tough" in Theorem 30 cannot be relaxed. On the other hand, the condition "$G$ is spider graph" in Theorem 30 cannot be removed since there are 3/2-tough nonhamiltonian graphs.\\

\subsection{Dominating cycles}

\noindent\textbf{Theorem 31}  (Nikoghosyan, 2011) \cite{[59]}

\noindent Let $G$ be a graph.  Then each longest cycle in $G$ is a dominating cycle if $\kappa\ge2$ and 
$$
q\le\left\{ 
\begin{array}{lll}
8 & \mbox{if} & \mbox{ }\delta=2, \\ \frac{3(\delta-1)(\delta+2)-1}{2} & \mbox{if} & \mbox{ }%
\delta\ge3. 
\end{array}
\right. 
$$

Examples for sharpness.   To show that Theorem 31 is sharp, suppose first
that $\delta = 2$. The graph $K_1 + 2K_2$ shows that the connectivity condition $\kappa\ge2$
in Theorem 31 cannot be relaxed by replacing it with $\kappa\ge 1$. The graph with
vertex set $\{v_1, v_2, ..., v_8\}$ and edge set
$$
\{v_1v_2, v_2v_3, v_3v_4, v_4v_5, v_5v_6, v_6v_1, v_1v_7, v_7v_8, v_8v_4\},
$$
shows that the size bound $q \le 8$ cannot be relaxed by replacing it with $q \le 9$.
Finally, the graph $K_2 +3K_1$ shows that the conclusion "each longest cycle in $G$
is a dominating cycle" cannot be strengthened by replacing it with "$G$ is hamiltonian".
Analogously, we can use $K_1 + 2K_{\delta}$, $K_2 + 3K_{\delta-1}$ and $K_{\delta} + (\delta + 1)K_1$,
respectively, to show that Theorem 31 is sharp when $\delta\ge 3$. So, Theorem 31 is
best possible in all respects.  \\

\noindent\textbf{Theorem 32} (Nash-Williams, 1971) \cite{[52]}
 
\noindent Let $G$ be a graph.   Then each longest cycle in $G$ is a dominating cycle if $\kappa\ge 2$ and
$$
\delta\ge \frac{n+2}{3}.
$$

Examples for sharpness. $2K_3+K_1$; $3K_{\delta-1}+K_2$; $H(1,2,4,3)$.\\

The graph $2K_3+K_1$ shows that the connectivity condition $\kappa\geq 2$ in Theorem 32 cannot be replaced by $\kappa\ge1$. The second graph shows that the minimum degree condition $\delta\geq (n+2)/3$ cannot be replaced by $\delta\geq (n+1)/2$. Finally, the third graph shows that the conclusion "is a dominating cycle" cannot be strengthened by replacing it with "is a Hamilton cycle". \\

\noindent\textbf{Theorem 33} (Bigalke and Jung, 1979) \cite{[9]}

\noindent Let $G$ be a graph.   Then  each longest cycle in $G$ is a dominating cycle if $\tau\ge 1$ and
$$
\delta\ge\frac{n}{3}.
$$

Examples for sharpness. $2(\kappa+1)K_2+\kappa K_1$; $L_3$; $G^*_n$.\\

\noindent\textbf{Theorem 34} (Yamashita, 2008) \cite{[69]} 

\noindent Let $G$ be graph.    Then each longest cycle in $G$ is a dominating cycle if $\kappa\ge 3$ and
$$
\delta\ge\frac{n+\kappa+3}{4}.
$$ 

Examples for sharpness. $3K_{\delta-1}+K_2$; $H(2,n-3\delta+3,\delta-1,\kappa)$; $H(1,2,\kappa+1,\kappa)$.\\

\subsection{$CD_{\lambda}$-cycles}

\noindent\textbf{Theorem 35} (Jung, 1990) \cite{[41]}

\noindent Let $G$ be a graph.   Then each longest cycle in $G$ is a $CD_3$-cycle if $\kappa\ge 3$ and
$$
\delta\ge \frac{n+6}{4}. 
$$

Examples for sharpness. $\lambda K_{\lambda+1}+K_{\lambda-1}$ $(\lambda\geq 2)$ ; $(\lambda+1)K_{\delta-\lambda+1}+K_{\lambda}$ $(\lambda\geq 1)$ ; $H(\lambda-1,\lambda,\lambda+2,\lambda+1)$ $(\lambda\geq 2)$. \\

\noindent\textbf{Theorem 36} (Nikoghosyan, 2009) \cite{[57]} 

\noindent Let $G$ be a graph and  $\lambda$ a positive integer.    Then each longest cycle in $G$ is a  $CD_{\min\{\lambda,\delta-\lambda+1\}}$-cycle if $\kappa\ge \lambda$ and 
$$
\delta\ge\frac{n+2}{\lambda+1}+\lambda-2.  
$$

Examples for sharpness. $\lambda K_{\lambda+1}+K_{\lambda-1}$ $(\lambda\geq 2)$ ; $(\lambda+1)K_{\delta-\lambda+1}+K_{\lambda}$ $(\lambda\geq 1)$ ; $H(\lambda-1,\lambda,\lambda+2,\lambda+1)$ $(\lambda\geq 2)$. \\

\subsection{Long cycles}

\noindent\textbf{Theorem 37} (Dirac, 1952) \cite{[23]} 

\noindent In every graph,  
$$
c\ge \delta+1.
$$

Example for sharpness. Join two copies of $K_{\delta+1}$ by an edge.\\

\noindent\textbf{Theorem 38} (Kouider, 1994) \cite{[47]} 

\noindent In every graph,  
$$
c\ge \frac{n}{\left\lceil \alpha/ \kappa\right\rceil}.
$$

Example for sharpness. Complete bipartite graph with $\kappa=\alpha$ shows that the bound in Theorem 38 is sharp. The original result is formulated for 2-connected graphs. However, Theorem 38 is true under assumption that each vertex (edge) is a cycle of length one (two, respectively). \\

\noindent\textbf{Theorem 39} (Nikoghosyan, 1998) \cite{[61]} 

\noindent Let $G$ be a graph and $C$ a longest cycle in $G$. Then 
$$
|C|\ge(\overline{p}+2)(\delta-\overline{p}).
$$  

Example for sharpness. $(\kappa+1)K_{\delta-\kappa+1}+K_\kappa$.\\

\noindent\textbf{Theorem 40} (Nikoghosyan, 2000) \cite{[61]} 

\noindent Let $G$ be a graph and $C$ a longest cycle in $G$. Then 
$$
|C|\ge(\overline{c}+1)(\delta-\overline{c}+1).
$$   

Example for sharpness. $(\kappa+1)K_{\delta-\kappa+1}+K_\kappa$.\\

\noindent\textbf{Theorem 41} (Nikoghosyan, 2000) \cite{[56]} 

\noindent Let $G$ be a graph with $\kappa\ge2$ and $C$ a longest cycle in $G$. If $\overline{c}\ge \kappa$ then 
$$
|C|\ge \frac{(\overline{c}+1)\kappa}{\overline{c}+\kappa+1}(\delta+2).
$$
Otherwise, 
$$
|C|\ge\frac{(\overline{c}+1)\overline{c}}{2\overline{c}+1}(\delta+2).
$$

Example for sharpness. $(\kappa+1)K_{\delta-\kappa+1}+K_\kappa$.\\

\subsection{Hamilton cycles and long cycles}

\noindent\textbf{Theorem 42} (Woodall, 1976) \cite{[67]} 

\noindent Let $G$ be a graph and $\lambda, t, r$ be integers with $n=t(\lambda-1)+r+1$, where $\lambda\ge2$, $t\ge0$ and $0\le r<\lambda-1$. If 
$$
q>t\left(^\lambda_2\right)+\left(^{r+1}_2\right)
$$
then
$$
c>\lambda.
$$

Example for sharpness. The result is best possible, in view of the graph consisting of $t$ copies of $K_\lambda$ and one copy of $K_{r+1}$, all having exactly one vertex in common.\\

\noindent\textbf{Theorem 43} (Fan, Lv and Wang, 2004) \cite{[30]}

\noindent Let $G$ be a 2-connected graph and let $2\le \lambda\le n-1$. If
$$
q>\max\left\{f(n,2,\lambda), f(n, \left\lfloor\frac{\lambda}{2}\right\rfloor, \lambda)\right\}
$$
then
$$
c>\lambda,
$$
where $f(n,t,\lambda)=(\lambda+1-t)(\lambda-t)/2+t(n-\lambda-1+t)$ and $2\le t\le \lambda/2$.\\

Examples for sharpness. The result is best possible, in view of the graph obtained from $K_{\lambda+1-t}$ by adding $n-(\lambda+1-t)$ isolated vertices, each joined to the same $t$ vertices of $K_{\lambda+1-t}$.\\

\noindent\textbf{Theorem 44} (Alon, 1986) \cite{[1]} 

\noindent Let $G$ be a graph and $\lambda$  a positive integer. If $\delta\ge\frac{n}{\lambda+1}$ then  
$$
c\ge\frac{n}{\lambda}.
$$

Examples for sharpness. $(\lambda+1)K_\lambda+K_1$; $\lambda K_{\lambda+1}$.\\

\noindent\textbf{Theorem 45} (Dirac, 1952) \cite{[23]}

\noindent Let $G$ be a graph. If $\kappa\ge 2$ then 
$$
c\ge\min\{n, 2\delta\}.
$$

Examples for sharpness. $(\lambda+1)K_{\lambda+1}+K_{\lambda}$ $(\lambda\geq 1)$; $(\lambda+3)K_{\lambda-1}+K_{\lambda+2}$ $(\lambda\geq 2)$; $(\lambda+2)K_{\lambda}+K_{\lambda+1}$ $(\lambda\geq 1)$.\\

\noindent\textbf{Theorem 46} (Kaneko and Yoshimoto, 1952) \cite{[44]}

\noindent Let $G$ be a 2-connected balanced bipartite graph. Then
$$
c\ge\min\{n, 4\delta-2\}.
$$

Examples for sharpness. Clearly, the condition "$G$ is balanced" in Theorem 46 cannot be removed. Consider the balanced bipartite graph $G=(X,Y;E)$ with vertex classes of the form $X=P\cup Q$, $Y=R\cup S$ with $z\in Q$, where $|P|=|R|=|Q|=|S|=n/4$, $N_G(x)=R$ for all $x\in P$, $N_G(x)=S$ for all $x\in Q-z$ and $N_G(z)=Y$. This example shows that 2-connectivity condition in Theorem 46 cannot be weakened. Next, consider the balanced bipartite graph $G=(X,Y;E)$ with vertex classes of the form $X=P\cup Q$, $Y=R\cup S$, where $|P|=|R|=|Q|=|S|=n/4$, $N_G(x)=R$ for all $x\in P$, and $N_G(x)=Y$ for all $x\in Q$. This example shows that the bound $4\delta-2$ in Theorem 46 cannot be improved.\\

\noindent\textbf{Theorem 47} (Bauer and Schmeichel, 1987) \cite{[4]}

\noindent Let $G$ be a graph. If $\tau\ge 1$ then 
$$
c\ge\min\{n, 2\delta+2\}.
$$

Examples for sharpness. $K_{\delta,\delta+1}$; $L_2$.\\

\noindent\textbf{Theorem 48} (Nikoghosyan, 2012) \cite{[60]}

\noindent Let $G$ be a graph. If $\tau > 4/3$ then 
$$
c\ge\min\{n, 2\delta+5\}.
$$

Examples for sharpness. The Petersen graph shows that the condition $\tau>4/3$ in Theorem 48 cannot be replaced by $\tau=4/3$. Let $H_1$ be a complete bipartite graph with bipartition $V_1=\{x_1,x_2,x_3,x_4,x_5\}$ and $V_2=\{y_1,y_2\}$, and let $H_2$ be a complete graph with vertex set $V=\{z_1,z_2,z_3,z_4,z_5\}$. The graph obtained from disjoint graphs $H_1$ and $H_2$ by adding the edges $x_iz_i$ $(i=1,...,5)$, shows that the bound $c\ge 2\delta+5$ in Theorem 48 cannot be replaced by $c\ge 2\delta+6$.\\

\noindent\textbf{Theorem 49} (Nikoghosyan, 1981) \cite{[54]} 

\noindent Let $G$ be a graph. If  $\kappa\ge3$ then  
$$
c\ge\min\{n,3\delta-\kappa\}.
$$

Examples for sharpness. $3K_{\delta-1}+K_2$; $H(1,\delta-\kappa+1,\delta,\kappa)$.\\

\noindent\textbf{Theorem 50} (Jung, 1978) \cite{[39]}

\noindent Let $G$ be a graph. If $\kappa\ge 3$ and  $\delta\ge \alpha$ then 
$$
c\ge \min\{n, 3\delta-3\}.
$$

Examples for sharpness. $(\lambda +2)K_{\lambda +2}+K_{\lambda +1}$; $(\lambda +4)K_{\lambda}+K_{\lambda +3}$; $(\lambda +3)K_{\lambda+1}+K_{\lambda +2}$.\\

\noindent\textbf{Theorem 51} (Nikoghosyan, 2009) \cite{[57]} 

\noindent Let $G$ be a graph and $\lambda$  a positive integer. If $\kappa\ge \lambda+2$  and  $\delta\ge \alpha+\lambda-1$ then  
$$
c\ge\min\{n,(\lambda+2)(\delta-\lambda)\}.
$$

Examples for sharpness. $(\lambda +2)K_{\lambda +2}+K_{\lambda +1}$; $(\lambda +4)K_{\lambda}+K_{\lambda +3}$; $(\lambda +3)K_{\lambda+1}+K_{\lambda +2}$.\\

\noindent\textbf{Theorem 52} (M.Zh. Nikoghosyan and Zh.G. Nikoghosyan, 2011) \cite{[53]} 

\noindent Let $G$ be a graph. If $\kappa\ge 4$  and  $\delta\ge \alpha$ then  
$$
c\ge\min\{n,4\delta-\kappa-4\}.
$$

Examples for sharpness. $4K_{\delta-2}+K_3$; $H(1,2,\kappa+1,\kappa)$; $H(2,n-3\delta+3,\delta-1,\kappa)$.\\

\noindent\textbf{Theorem 53} (Bauer, Morgana, Schmeichel and Veldman, 1989) \cite{[3]}

\noindent Let $G$ be a graph. If $\kappa\ge 2$ and  $\delta\ge \frac{n+2}{3}$ then 
$$
c\ge \min\{n, n+\delta-\alpha\}.
$$

Examples for sharpness. $2K_\delta+K_1$; $3K_{\delta-1}+K_2$; $K_{2\delta-2,\delta}$.\\

\noindent\textbf{Theorem 54} (Bauer, Schmeichel and Veldman, 1988) \cite{[6]}

\noindent Let $G$ be a graph. If $\tau\ge 1$ and $\delta\ge \frac{n}{3}$ then 
$$
c\ge \min\{n, n+\delta-\alpha+1\}.
$$

Examples for sharpness. $K_{\delta,\delta+1}$; $L_\delta$; $G^*_n$.\\

\subsection{Dominating cycles and long cycles}

\noindent\textbf{Theorem 55} (Jung, 1981) \cite{[40]} 

\noindent Let $G$ be a graph. If  $\kappa\ge 3$ then either  each longest cycle in $G$ is a dominating cycle or
$$
c\ge 3\delta-3.
$$
 
Examples for sharpness. $(\lambda+1)K_{\lambda+1}+K_{\lambda}$ $(\lambda\geq 1)$; $(\lambda+3)K_{\lambda-1}+K_{\lambda+2}$ $(\lambda\geq 2)$; $(\lambda+2)K_{\lambda}+K_{\lambda+1}$ $(\lambda\geq 1)$.\\

\noindent\textbf{Theorem 56} (M.Zh. Nikoghosyan and Zh.G. Nikoghosyan, 2011) \cite{[53]} 

\noindent Let $G$ be a graph. If  $\kappa\ge 4$ then either each longest cycle in $G$ is a dominating cycle or
$$
c\ge 4\delta-\kappa-4. 
$$             

Examples for sharpness. $4K_{\delta-2}+K_3$; $H(2,\delta-\kappa+1,\delta-1,\kappa)$; $H(1,2,\kappa+1,\kappa)$.

\subsection{$CD_{\lambda}$-cycles and long cycles}

\noindent\textbf{Theorem 57} (Nikoghosyan, 2009) \cite{[57]} 

\noindent Let $G$ be a graph and $\lambda$ a positive integer. If  $\kappa\ge \lambda+1$ then either each longest cycle in $G$ is a  $CD_{\min\{\lambda,\delta-\lambda\}}$-cycle or 
$$
c\ge (\lambda+1)(\delta-\lambda+1). 
$$    

Examples for sharpness. $(\lambda+1)K_{\lambda+1}+K_{\lambda}$ $(\lambda\geq 1)$; $(\lambda+3)K_{\lambda-1}+K_{\lambda+2}$ $(\lambda\geq 2)$; $(\lambda+2)K_{\lambda}+K_{\lambda+1}$ $(\lambda\geq 1)$.\\

\noindent Institute for Informatics and Automation Problems\\ National Academy of Sciences\\
P. Sevak 1, Yerevan 0014, Armenia\\ 
E-mail: zhora@ipia.sci.am
\end{document}